# Space-Time-Unified GIRM (Generalized Integral Representation Method) for Unsteady Advective Diffusion

Hiroshi Isshiki, Institute of Mathematical Analysis, Osaka, Japan,
isshiki@dab.hi-ho.ne.jp

Daisuke Kitazawa, Institute of Industrial Science,
The University of Tokyo, dkita@iis.u-tokyo.ac.jp

**Abstract**
The Generalized Integral Representation Method (GIRM) for Space-Time-Separated Method (STSM) and Space-Time-Unified Method (STUM) are discussed. STSM and STUM give explicit and implicit time evolutions, respectively. The algorithm of STSM is much simpler than STUM. However, the implicit time evolution of STUM could give us much more efficient computation. Numerical calculations using STUM for Dirichlet and Neumann problems in 2D space-time are conducted using a Traditional Fundamental Solution (TFS). The results seem very satisfactory. However, in case of Neumann problem, the time increment must be smaller than in case of Dirichlet problem.

**1. Introduction**

There are two kinds of approaches to the time evolution of the initial-boundary value problems:
 (1) Space-Time-Separated Methods (STSM)
 (2) Space-Time-Unified Methods (STUM)
In STSM, space and time treated separately. Concretely, the present time step predicts the time derivatives. The next time step is estimated explicitly using the time derivatives. In STUM, the next time step is estimated implicitly using the present time step.

Historically, differential equations are used to obtain mathematical models. In order to obtain numerical solutions, we must discretize the differential equations. Finite Difference Method (FDM) is a strong solution that satisfies the differential equation point by point. The accuracy is very high. However, FDM can't be applied to irregular mesh in general. If we wish to use irregular mesh, we must introduce the concept of weak solution such as Finite Element Method (FEM) that satisfies the differential equation averagely. Generally speaking, the weak solution is very flexible, but the accuracy is lower. Namely, the computational time is longer.

On the other hand, if the differential expression is transformed into an integral expression, we can use an irregular mesh. Integral Representation Method (IRM) could use the irregular mesh and realize more flexible numerical method without sacrificing the accuracy. But, IRM is limited to the system of linear diff. equation with constant coefficients, since the Fundamental Solution (FS) is an analytical singular solution of the differential equation[1,2].

IRM is extended by the author to Generalized Integral Representation Method (GIRM)[3,4,5]. Generalized Fundamental Solution (GFS) is used in Generalized Integral Representation Method (GIRM). GFS is not required to satisfy the differential equation. Usually, element based methods require overall $C^0$ continuity. GIRM does not require the inter-element continuity. This could make the mesh generation very flexible.

**2. Steady Diffusion Problem**

**2D Steady Diffusion Problem**

For simplicity, we consider a 2D steady problem in a region $\Omega$ with the boundary $\Gamma$ as shown in Fig. 1. Let the coordinates be $(x_1, x_2)$ and $C(x_1, x_2)$ refer to the concentration of a substance, temperature or something like that. A steady diffusion problem is given in an ordinary notation as

$$U_1(x_1,x_2)\frac{\partial C(x_1,x_2)}{\partial x_1}+U_2(x_1,x_2)\frac{\partial C(x_1,x_2)}{\partial x_1}$$
$$=\frac{\partial}{\partial x_1}\left(\nu(x_1,x_2)\frac{\partial C(x_1,x_2)}{\partial x_1}\right)+\frac{\partial}{\partial x_2}\left(\nu(x_1,x_2)\frac{\partial C(x_1,x_2)}{\partial x_2}\right)+f(x_1,x_2), \quad (x_1,x_2)\in\Omega \quad (1)$$

$$\frac{\partial C(x_1,x_2)}{\partial n}+a(x_1,x_2)C(x_1,x_2)=b(x_1,x_2), \quad (x_1,x_2)\in\Gamma \quad (2)$$

where $(U_1(x_1,x_2), U_2(x_1,x_2))$ and $\nu(x_1,x_2)$ are the advective velocity and diffusion constant, respectively. Functions $f(x_1,x_2)$, $a(x_1,x_2)$ and $b(x_1,x_2)$ are given functions of the coordinates $(x_1,x_2)$.

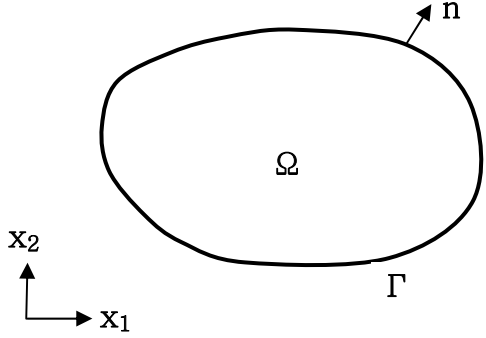

Fig. 1. 2D region and boundary

In this paper, we use the following notation:

$$U_i\frac{\partial C(x)}{\partial x_i}=\frac{\partial}{\partial x_i}\left(\frac{\partial C(x)}{\partial x_i}\right)+\sigma(x), \quad x\in\Omega \quad (3)$$

$$\frac{\partial C(x)}{\partial n(x)}+a(x)C(x)=b(x), \quad x\in\Gamma \quad (4)$$

We introduce a Generalized Fundamental Solution (GFS) $G(x,\xi)$. In case of a conventional Fundamental Solution (FS), we require $G(x,\xi)$ as a singular fundamental solution of a differential equation, for example (see Appendix A):

$$\Delta G(x,\xi)=\boldsymbol{\delta}(x-\xi)=\delta(x_1-\xi_1)\delta(x_2-\xi_2) \quad (5)$$

$$G(x,\xi)=\frac{1}{2\pi}\ln\sqrt{(x_1-\xi_1)^2+(x_2-\xi_2)^2} \quad (6)$$

However, in case of GFS, we require only $C^{(2)}$ continuity.

Firstly, we obtain an Integral Representation of the basic equations (3) and (4) using GFS $G(x,\xi)$:

$$0 = \int_{\Omega(\xi)} \left( U_i(\xi) \frac{\partial C(\xi)}{\partial \xi_i} - \frac{\partial}{\partial \xi_i}\left(\nu(\xi)\frac{\partial C(\xi)}{\partial \xi_i}\right) + \sigma(\xi) \right) G(\xi, x) d\Omega(\xi)$$

$$= \int_{\Omega(\xi)} \left( \frac{\partial U_i(\xi)C(\xi)G(\xi,x)}{\partial \xi_i} - C(\xi)\frac{\partial U_i(\xi)G(\xi,x)}{\partial \xi_i} - \frac{\partial}{\partial \xi_i}\left(\nu(\xi)\frac{\partial C(\xi)}{\partial \xi_i}G(\xi,x)\right) + \nu(\xi)\frac{\partial C(\xi)}{\partial \xi_i}\frac{\partial G(\xi,x)}{\partial \xi_i} - \sigma(\xi)G(\xi,x) \right) d\Omega(\xi)$$

$$= \int_{\Gamma(\xi)} U_i(\xi)n_i(\xi)C(\xi)G(\xi,x) d\Gamma(\xi) - \int_{\Omega(\xi)} C(\xi)\frac{\partial U_i(\xi)G(\xi,x)}{\partial \xi_i} d\Omega(\xi) - \int_{\Gamma(\xi)} \nu(\xi)\frac{\partial C(\xi)}{\partial n(\xi)}G(\xi,x) d\Gamma(\xi) ,\quad (7)$$

$$+ \int_{\Omega(\xi)} \left( \frac{\partial}{\partial \xi_i}\left(C(\xi)\nu(\xi)\frac{\partial G(\xi,x)}{\partial \xi_i}\right) - C(\xi)\frac{\partial}{\partial \xi_i}\left(\nu(\xi)\frac{\partial G(\xi,x)}{\partial \xi_i}\right) - \sigma(\xi)G(\xi,x) \right) d\Omega(\xi)$$

$$= \int_{\Gamma(\xi)} \left( -\left(\nu(\xi)\frac{\partial C(\xi)}{\partial n(\xi)} - U_i(\xi)n_i(\xi)C(\xi)\right)G(\xi,x) + C(\xi)\nu(\xi)\frac{\partial G(\xi,x)}{\partial n(\xi)} \right) d\Gamma(\xi)$$

$$- \int_{\Omega(\xi)} \left( C(\xi)\frac{\partial U_i(\xi)G(\xi,x)}{\partial \xi_i} + C(\xi)\frac{\partial}{\partial \xi_i}\left(\nu(\xi)\frac{\partial G(\xi,x)}{\partial \xi_i}\right) + \sigma(\xi)G(\xi,x) \right) d\Omega(\xi)$$

where $n(\xi)$ is the unit outward normal of $\Gamma(\xi)$. Substituting Eq. (4) into Eq. (7), we have

$$0 = \int_{\Gamma(\xi)} \left( -\left(\nu(\xi)(b(\xi) - a(\xi)C(\xi)) - U_i(\xi)n_i(\xi)C(\xi)\right)G(\xi,x) + C(\xi)\nu(\xi)\frac{\partial G(\xi,x)}{\partial n(\xi)} \right) d\Gamma(\xi)$$
$$- \int_{\Omega(\xi)} \left( C(\xi)\frac{\partial U_i(\xi)G(\xi,x)}{\partial \xi_i} + C(\xi)\frac{\partial}{\partial \xi_i}\left(\nu(\xi)\frac{\partial G(\xi,x)}{\partial \xi_i}\right) + \sigma(\xi)G(\xi,x) \right) d\Omega(\xi) \quad (8)$$

Eq. (8) is a Fredholm Integral Equation of the first kind with unknown $C(x)$. Namely, $C(x)$, $x \in \Omega(x) + \Gamma(x)$ can be determined solving Eq. (8).

If the diffusion coefficient $\nu(\xi)$ is constant, we can use the conventional FS defined by Eq. (5). Substituting Eq. (6) into Eq. (8), we obtain

$$\varepsilon(x)C(x) = \int_{\Gamma(\xi)} \left( -\left((b(\xi) - a(\xi)C(\xi)) - U_i(\xi)n_i(\xi)C(\xi)\right)G(\xi,x) + C(\xi)\frac{\partial G(\xi,x)}{\partial n(\xi)} \right) d\Gamma(\xi)$$
$$- \frac{1}{\nu} \int_{\Omega(y)} \left( C(\xi)\frac{\partial U_i(\xi)G(\xi,x)}{\partial \xi_i} + \sigma(\xi)G(\xi,x) \right) d\Omega(\xi) \quad (9)$$

where

$$\varepsilon(x) = 1, \text{ if } x \in \Omega; \quad 1/2, \text{ if } x \in \Gamma; \quad 0, \text{ otherwise} \quad (10)$$

When $x \in \Gamma(x)$ and $U_i(x)$ is zero, Eq. (9) becomes

$$\frac{1}{2}C(x) = \int_{\Gamma(\xi)} \left( -\left((b(\xi) - a(\xi)C(\xi))\right)G(\xi,x) + C(\xi)\frac{\partial G(\xi,x)}{\partial n(\xi)} \right) d\Gamma(\xi) - \frac{1}{\nu}\int_{\Omega(y)} \sigma(\xi)G(\xi,x) d\Omega(\xi) \quad (11)$$

Eq. (14) is a Fredholm Integral Equation of the Second kind with unknown $C(x)$ on $\Gamma(x)$. Namely, $C(x)$, $x \in \Gamma(x)$ can be determined solving Eq. (11). $C(x)$, $x \in \Omega(x)$ can be obtained by Eq. (9).

## 3. Unsteady Diffusion Problem

### 3.1. 1D Unsteady Diffusion Problem

For simplicity, we consider first a 1D (or 2D in space-time) unsteady problem in a region $(-L, +L)$ with the boundary $x = \pm L$ as shown in Fig. 2. Let the space-time coordinates be $(x, t)$ and $C(x, t)$ refer to the concentration of a substance. An unsteady advective diffusion problem is given as

$$\frac{\partial C(x,t)}{\partial t} + U\frac{\partial C(x,t)}{\partial x} = \frac{\partial}{\partial x}\left(\nu(x,t)\frac{\partial C(x,t)}{\partial x}\right) + \sigma(x,t), \quad -L < x < L, \ 0 < t < T \tag{1}$$

$$\frac{\partial C(x,t)}{\partial n(x)} + a(x,t)C(x,t) = b(x,t), \quad x = \pm L, \ 0 < t < T \tag{2}$$

$$C(x,0) = C_0(x), \quad -L < x < +L \tag{3}$$

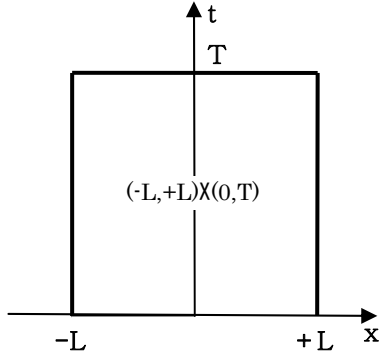

Fig. 2. 2D space-time

Firstly, we obtain an Integral Representation of the basic equations (1), (2) and (3)[6] using adjoint Generalized Fundamental Solution (GFS) $G*(x,\xi,t,\tau)$ from

$$0 = \int_0^T \int_{-L}^{+L}\left(\frac{\partial C(\xi,\tau)}{\partial \tau} + U(\xi,\tau)\frac{\partial C(\xi,\tau)}{\partial \xi} - \frac{\partial}{\partial \xi}\left(\nu(\xi,\tau)\frac{\partial C(\xi,\tau)}{\partial \xi}\right) - \sigma(\xi,\tau)\right)G*(\xi,x,\tau,t)d\xi d\tau \tag{4}$$

where GFS $G(x,\xi,t,\tau)$ and adjoint GFS $G*(x,\xi,t,\tau)$ satisfies the reverse causality:

$$G(x,\xi,t,\tau) = 0 \ \text{when} \ t < \tau \ \text{and} \ G*(x,\xi,t,\tau) = 0 \ \text{when} \ t > \tau \tag{5a}$$

$$G(x,\xi,t,\tau) = G*(\xi,x,\tau,t) \tag{5b}$$

In case of a Traditional Fundamental Solution (TFS) with $\nu = const$, we require $G(x,\xi,t,\tau)$ as a singular fundamental solution of a differential equation, for example (see Appendix B):

$$\frac{\partial G(x,\xi,t,\tau)}{\partial t} - \nu \frac{\partial^2 G(x,\xi,t,\tau)}{\partial x^2} = \delta(x-\xi)\delta(t-\tau), \quad (x,t) \in R^2 \times (0,\infty) \tag{6a}$$

$$G(x,\xi,t,\tau) = \frac{1}{\sqrt{4\pi\nu(t-\tau)}}\exp\left(-\frac{|x-\xi|^2}{4\nu(t-\tau)}\right)H(t-\tau) = G(\xi,x,t,\tau) = G(\xi,x,-\tau,-t) \tag{6b}$$

and

$$-\frac{\partial G*(x,\xi,t,\tau)}{\partial t} - \nu \frac{\partial^2 G*(x,\xi,t,\tau)}{\partial x^2} = \delta(x-\xi)\delta(t-\tau), \quad (x,t) \in R^2 \times (0,\infty) \tag{7a}$$

$$G*(x,\xi,t,\tau) = \frac{1}{\sqrt{4\pi\nu(\tau-t)}}\exp\left(-\frac{|\xi-x|^2}{4\nu(\tau-t)}\right)H(\tau-t) = G(\xi,x,\tau,t) \tag{7b}$$

We rewrite Eq. (4) using integration by part. The each term in Eq. (4) becomes as follows:

$$\int_0^T \int_{-L}^{+L}\frac{\partial C(\xi,\tau)}{\partial \tau}G*(\xi,x,\tau,t)d\xi d\tau$$

$$= \int_0^T \int_{-L}^{+L}\frac{\partial C(\xi,\tau)G*(\xi,x,\tau,t)}{\partial \tau}d\xi d\tau - \int_0^T \int_{-L}^{+L}C(\xi,\tau)\frac{\partial G*(\xi,x,\tau,t)}{\partial \tau}d\xi d\tau \tag{8a}$$

$$= -\int_{-L}^{+L}C(\xi,0)G*(\xi,x,0,t)d\xi - \int_0^T \int_{-L}^{+L}C(\xi,\tau)\frac{\partial G*(\xi,x,\tau,t)}{\partial \tau}d\xi d\tau$$

$$\int_0^T \int_{-L}^{+L}U(\xi,\tau)\frac{\partial C(\xi,\tau)}{\partial \xi}G*(\xi,x,\tau,t)d\xi d\tau$$

$$= \int_0^T \int_{-L}^{+L}\frac{\partial U(\xi,\tau)C(\xi,\tau)G*(\xi,x,\tau,t)}{\partial \xi}d\xi d\tau - \int_0^T \int_{-L}^{+L}C(\xi,\tau)\frac{\partial U(\xi,\tau)G*(\xi,x,\tau,t)}{\partial \xi}d\xi d\tau \tag{8b}$$

$$= \int_0^T \left[U(\xi,\tau)C(\xi,\tau)G*(\xi,x,\tau,t)\right]_{\xi=-L}^{\xi=+L}d\tau - \int_0^T \int_{-L}^{+L}C(\xi,\tau)\frac{\partial U(\xi,\tau)G*(\xi,x,\tau,t)}{\partial \xi}d\xi d\tau$$

$$-\int_0^T \int_{\Omega(\xi)} \frac{\partial}{\partial \xi}\left(\nu(\xi,\tau)\frac{\partial C(\xi,\tau)}{\partial \xi}\right) G*(\xi,x,\tau,t)\,d\Omega(\xi)\,d\tau$$

$$=-\int_0^T \int_{-L}^{+L} \frac{\partial}{\partial \xi}\left(\nu(\xi,\tau)\frac{\partial C(\xi,\tau)}{\partial \xi} G*(\xi,x,\tau,t)\right) d\xi\,d\tau + \int_0^T \int_{-L}^{+L} \nu(\xi,\tau)\frac{\partial C(\xi,\tau)}{\partial \xi}\frac{\partial G*(\xi,x,\tau,t)}{\partial \xi}\,d\xi\,d\tau$$

$$=-\int_0^T \left[\nu(\xi,\tau)\frac{\partial C(\xi,\tau)}{\partial \xi} G*(\xi,x,\tau,t)\right]_{\xi=-L}^{\xi=+L} d\tau + \int_0^T \int_{-L}^{+L} \frac{\partial}{\partial \xi}\left(C(\xi,\tau)\nu(\xi,\tau)\frac{\partial G*(\xi,x,\tau,t)}{\partial \xi}\right) d\xi\,d\tau$$

$$-\int_0^T \int_{-L}^{+L} C(\xi,\tau)\frac{\partial}{\partial \xi}\left(\nu(\xi,\tau)\frac{\partial G*(\xi,x,\tau,t)}{\partial \xi}\right) d\xi\,d\tau \quad (8c)$$

$$=-\int_0^T \left[\nu(\xi,\tau)\left(\frac{\partial C(\xi,\tau)}{\partial \xi} G*(\xi,x,\tau,t) - C(\xi,\tau)\frac{\partial G*(\xi,x,\tau,t)}{\partial \xi}\right)\right]_{\xi=-L}^{\xi=+L} d\tau$$

$$-\int_0^T \int_{-L}^{+L} C(\xi,\tau)\frac{\partial}{\partial \xi}\left(\nu(\xi,\tau)\frac{\partial G*(\xi,x,\tau,t)}{\partial \xi}\right) d\xi\,d\tau$$

Substituting Eq. (8) into Eq. (4), we have

$$0 = -\int_{-L}^{+L} C(\xi,0) G*(\xi,x,0,t)\,d\xi$$

$$+ \int_0^T \left[\begin{array}{l} U(\xi,\tau)C(\xi,\tau)G*(\xi,x,\tau,t) \\ -\nu(\xi,\tau)\left(\frac{\partial C(\xi,\tau)}{\partial \xi} G*(\xi,x,\tau,t) - C(\xi,\tau)\frac{\partial G*(\xi,x,\tau,t)}{\partial \xi}\right) \end{array}\right]_{\xi=-L}^{\xi=+L} d\tau \quad (9)$$

$$-\int_0^T \int_{-L}^{+L} \left[\begin{array}{l} C(\xi,\tau)\left(\frac{\partial G*(\xi,x,\tau,t)}{\partial \tau} + \frac{\partial U(\xi,\tau)G*(\xi,x,\tau,t)}{\partial \xi} + \frac{\partial}{\partial \xi}\left(\nu(\xi,\tau)\frac{\partial G*(\xi,x,\tau,t)}{\partial \xi}\right)\right) \\ + \sigma(\xi,\tau)G*(\xi,x,\tau,t) \end{array}\right] d\xi\,d\tau$$

Substituting Eqs. (2), (3) and (7b) into Eq. (9), we obtain

$$0 = -\int_{-L}^{+L} C_0(\xi) G(x,\xi,t,0)\,d\xi$$

$$+ \int_0^t \left[\begin{array}{l} (U(\xi,\tau)C(\xi,\tau) - \nu(\xi,\tau)[b(\xi,\tau) - a(\xi,\tau)C(\xi,\tau)])G(x,\xi,t,\tau) \\ + \nu(\xi,\tau)C(\xi,\tau)\frac{\partial G(x,\xi,t,\tau)}{\partial \xi} \end{array}\right]_{\xi=-L}^{\xi=+L} d\tau$$

$$-\int_0^t \int_{-L}^{+L} \left[\begin{array}{l} C(\xi,\tau)\left(\frac{\partial G(x,\xi,t,\tau)}{\partial \tau} + \frac{\partial U(\xi,\tau)G(x,\xi,t,\tau)}{\partial \xi} + \frac{\partial}{\partial \xi}\left(\nu(\xi,\tau)\frac{\partial G(x,\xi,t,\tau)}{\partial \xi}\right)\right) \\ + \sigma(\xi,\tau)G(x,\xi,t,\tau) \end{array}\right] d\xi\,d\tau$$

$$= -\int_{-L}^{+L} C_0(\xi) G(x,\xi,t,0)\,d\xi$$

$$+ \int_0^t \left[\begin{array}{l} -(U(-L,\tau)C(-L,\tau) - \nu(-L,\tau)[b(-L,\tau) - a(-L,\tau)C(-L,\tau)])G(x,-L,t,\tau) \\ + \nu(-L,\tau)C(-L,\tau)\frac{\partial G(x,-L,t,\tau)}{\partial x} \\ + (U(+L,\tau)C(+L,\tau) - \nu(+L,\tau)[b(+L,\tau) - a(+L,\tau)C(+L,\tau)])G(x,+L,t,\tau) \\ - \nu(+L,\tau)C(+L,\tau)\frac{\partial G(x,+L,t,\tau)}{\partial x} \end{array}\right] d\tau$$

$$-\int_0^t \int_{-L}^{+L} \left[\begin{array}{l} C(\xi,\tau)\left(\begin{array}{l} -\frac{\partial G(x,\xi,t,\tau)}{\partial t} + \frac{\partial U(\xi,\tau)}{\partial \xi} G(x,\xi,t,\tau) - U(\xi,\tau)\frac{\partial G(x,\xi,t,\tau)}{\partial x} \\ -\frac{\partial \nu(\xi,\tau)}{\partial \xi}\frac{\partial G(x,\xi,t,\tau)}{\partial x} + \nu(\xi,\tau)\frac{\partial^2 G(x,\xi,t,\tau)}{\partial x^2} \end{array}\right) \\ + \sigma(\xi,\tau)G(x,\xi,t,\tau) \end{array}\right] d\xi\,d\tau \quad , \quad (10)$$

where the integral limit $T$ is replaced by $t$ because of Heaviside function $H(t-\tau)$ in $G$, and . $\partial G/\partial \tau = -\partial G/\partial t$ and $\partial G/\partial \xi = -\partial G/\partial x$ are used. Eq. (10) is a Volterra-Fredholm Integral Equation with unknown $C(x,t)$. Namely, $C(x,t)$, $-L \le x \le +L$, $0 < t$ can be determined solving Eq. (10).

If the diffusion coefficient $\nu(x,t)$ is constant and TFS $G(x,\xi,t,\tau)$ given by Eq. (6) is used, we obtain

$$\begin{aligned}\varepsilon(x,t)C(x,t) &= \int_{-L}^{+L} C_0(\xi)G(x,\xi,t,0)\,d\xi \\ &- \int_0^t \begin{bmatrix} -(U(-L,\tau)C(-L,\tau)-\nu(b(-L,\tau)-a(-L,\tau)C(-L,\tau)))G(x,-L,t,\tau)+\nu C(-L,\tau)\dfrac{\partial G(x,-L,t,\tau)}{\partial x} \\ +(U(+L,\tau)C(+L,\tau)-\nu(b(+L,\tau)-a(+L,\tau)C(+L,\tau)))G(x,+L,t,\tau)-\nu C(+L,\tau)\dfrac{\partial G(x,+L,t,\tau)}{\partial x} \end{bmatrix} d\tau \\ &+ \int_0^t \int_{-L}^{+L} \begin{bmatrix} C(\xi,\tau)\left(+\dfrac{\partial U(\xi,\tau)}{\partial \xi}G(x,\xi,t,\tau)-U(\xi,\tau)\dfrac{\partial G(x,\xi,t,\tau)}{\partial x}\right) \\ +\sigma(\xi,\tau)G(x,\xi,t,\tau) \end{bmatrix} d\xi d\tau \end{aligned} \quad (11)$$

where

$$\varepsilon(x,t)=\begin{cases} 1 & \text{if } (x,t)\in\{(x,t)\,|\,-L<x<+L,\,0<t<T\} \\ 1/2 & \text{if } (x,t)\in\{(x,t)\,|\,-L<x<+L,\,t=0\} \\ 1/2 & \text{if } (x,t)\in\{(x,t)\,|\,x=\pm L,\,0<t<T\} \\ 0 & \text{otherwise} \end{cases} \quad (12)$$

Furthermore, if $U(x)=0$, we have

$$\begin{aligned}\varepsilon(x,t)C(x,t) &= \int_{-L}^{+L} C_0(\xi)G(x,\xi,t,0)\,d\xi \\ &- \int_0^t \begin{bmatrix} +\nu[b(-L,\tau)-a(-L,\tau)C(-L,\tau)]G(x,-L,t,\tau)+\nu C(-L,\tau)\dfrac{\partial G(x,-L,t,\tau)}{\partial x} \\ -\nu[b(+L,\tau)-a(+L,\tau)C(+L,\tau)]G(x,+L,t,\tau)-\nu C(+L,\tau)\dfrac{\partial G(x,+L,t,\tau)}{\partial x} \end{bmatrix} d\tau \\ &+ \int_0^t \int_{-L}^{+L} \sigma(\xi,\tau)G(x,\xi,t,\tau)\,d\xi d\tau \end{aligned} \quad (13)$$

When $x=\pm L$, Eq. (13) becomes

$$\begin{aligned}\tfrac{1}{2}C(x,t) &= \int_{-L}^{+L} C_0(\xi)G(x,\xi,t,0)\,d\xi \\ &- \int_0^t \begin{bmatrix} +\nu[b(-L,\tau)-a(-L,\tau)C(-L,\tau)]G(x,-L,t,\tau)+\nu C(-L,\tau)\dfrac{\partial G(x,-L,t,\tau)}{\partial x} \\ -\nu[b(+L,\tau)-a(+L,\tau)C(+L,\tau)]G(x,+L,t,\tau)-\nu C(+L,\tau)\dfrac{\partial G(x,+L,t,\tau)}{\partial x} \end{bmatrix} d\tau \\ &+ \int_0^t \int_{-L}^{+L} \sigma(\xi,\tau)G(x,\xi,t,\tau)\,d\xi d\tau \end{aligned} \quad (14)$$

Eq. (14) is a Volterra-Fredholm Integral Equation with unknown $C(x,t)$, $x=\pm L$, $t>0$. Namely, $C(x,t)$, $x=\pm L$, $t>0$ can be determined solving Eq. (14). $C(x,t)$, $-L<x<+L$, $t>0$ can be obtained by Eq. (13).

### 3.2. 2D or 3D Space-Time Unsteady Diffusion Problem

Now, we consider a nD (or (n+1)D in space-time) unsteady problem in a region $\Omega$ with the boundary $\Gamma$ as shown in Fig. 3 for 3D space-time. Let the space-time coordinates be $(x_1,\cdots,x_n,t)$ and $C(x_1,\cdots,x_n,t)$ refer to the concentration of a substance. An unsteady advective diffusion problem is given as

$$\frac{\partial C(x,t)}{\partial t}+U_i\frac{\partial C(x,t)}{\partial x_i}=\frac{\partial}{\partial x_i}\left(\frac{\partial C(x,t)}{\partial x_i}\right)+\sigma(x,t),\quad (x,t)\in\Omega\times(0,T) \quad (15)$$

$$\frac{\partial C(x,t)}{\partial n(x)}+a(x,t)C(x,t)=b(x,t),\quad (x,t)\in\Gamma\times(0,T) \quad (16)$$

$$C(x,0)=C_0(x),\quad x\in\Omega \quad (17)$$

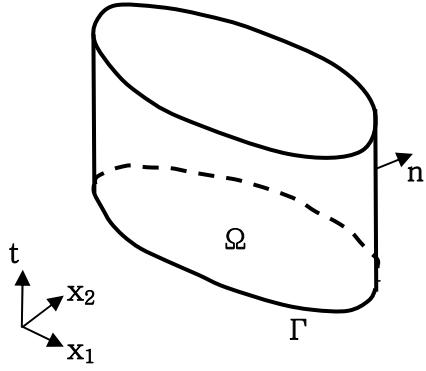

Fig. 3. 3D space-time

Firstly, we obtain an Integral Representation of the basic equations (12), (13) and (14)[6] using adjoint GFS $G*(x,\xi,t,\tau)$:

$$0 = \int_0^T \int_{\Omega(\xi)} \left( \frac{\partial C(\xi,\tau)}{\partial \tau} + U_i(\xi,\tau)\frac{\partial C(\xi,\tau)}{\partial \xi_i} - \frac{\partial}{\partial \xi_i}\left(\nu(\xi,\tau)\frac{\partial C(\xi,\tau)}{\partial \xi_i}\right) - \sigma(\xi,\tau) \right) G*(\xi,x,\tau,t)d\Omega(\xi)d\tau \tag{18}$$

where GFS $G(x,\xi,t,\tau)$ and adjoint GFS $G*(x,\xi,t,\tau)$ satisfies the reverse causality:

$$G(x,\xi,t,\tau) = 0 \text{ when } \tau > t \text{ and } G*(x,\xi,t,\tau) = 0 \text{ when } t > \tau \tag{19a}$$

$$G(x,\xi,t,\tau) = G*(\xi,x,\tau,t) \tag{19b}$$

In case of a Traditional Fundamental Solution (TFS) with $U_i = 0$, we require $G(x,\xi)$ as a singular fundamental solution of a differential equation, for example:

$$\frac{\partial G(x,\xi,t,\tau)}{\partial t} - \nu\Delta G(x,\xi,t,\tau) = \boldsymbol{\delta}(x-\xi)\delta(t-\tau) = \delta(x_1-\xi_1)\cdots\delta(x_n-\xi_n)\delta(t-\tau),\ (x,t)\in R^2\times(0,\infty) \tag{20a}$$

$$G(x,\xi,t,\tau) = \frac{1}{(4\pi\nu(t-\tau))^{n/2}}\exp\left(-|x-\xi|^2/(4\nu(t-\tau))\right)H(t-\tau) = G*(\xi,x,\tau,t) \tag{20b}$$

and

$$-\frac{\partial G*(x,\xi,t,\tau)}{\partial t} - \nu\Delta G*(x,\xi,t,\tau) = \boldsymbol{\delta}(x-\xi)\delta(t-\tau) = \delta(x_1-\xi_1)\cdots\delta(x_n-\xi_n)\delta(t-\tau),\ (x,t)\in R^2\times(0,\infty) \tag{21a}$$

$$G*(x,\xi,t,\tau) = \frac{1}{(4\pi\nu(t-\tau))^{n/2}}\exp\left(-\frac{|\xi-x|^2}{4\nu(\tau-t)}\right)H(\tau-t) = G(\xi,x,\tau,t) \tag{21b}$$

where $n$ is the number of space dimension, and $H(t)$ is Heaviside step function.

We rewrite Eq. (18) using integration by part. The each term in Eq. (18) becomes as follows:

$$\int_0^T \int_{\Omega(\xi)} \frac{\partial C(\xi,\tau)}{\partial \tau} G*(\xi,x,\tau,t)d\Omega(\xi)d\tau$$

$$= \int_0^T \int_{\Omega(\xi)} \frac{\partial C(\xi,\tau)G*(\xi,x,\tau,t)}{\partial \tau} d\Omega(\xi)d\tau - \int_0^T \int_{\Omega(\xi)} C(\xi,\tau)\frac{\partial G*(\xi,x,\tau,t)}{\partial \tau} d\Omega(\xi)d\tau \tag{22a}$$

$$= -\int_{\Omega(\xi)} C(\xi,0)G*(\xi,x,0,t)d\Omega(\xi) - \int_0^T \int_{\Omega(\xi)} C(\xi,\tau)\frac{\partial G*(\xi,x,\tau,t)}{\partial \tau} d\Omega(\xi)d\tau$$

$$\int_0^T \int_{\Omega(\xi)} U_i(\xi,\tau)\frac{\partial C(\xi,\tau)}{\partial \xi_i} G*(\xi,x,\tau,t)d\Omega(\xi)d\tau$$

$$= \int_0^T \int_{\Omega(\xi)} \frac{\partial U_i(\xi,\tau)C(\xi,\tau)G*(\xi,x,\tau,t)}{\partial \xi_i} d\Omega(\xi)d\tau - \int_0^T \int_{\Omega(\xi)} C(\xi,\tau)\frac{\partial U_i(\xi,\tau)G*(\xi,x,\tau,t)}{\partial \xi_i} d\Omega(\xi)d\tau \tag{22b}$$

$$= \int_0^T \int_{\Gamma(\xi)} U_i(\xi,\tau)n_i(\xi,\tau)C(\xi,\tau)G*(\xi,x,\tau,t)d\Gamma(\xi)d\tau - \int_0^T \int_{\Omega(\xi)} C(\xi,\tau)\frac{\partial U_i(\xi,\tau)G*(\xi,x,\tau,t)}{\partial \xi_i} d\Omega(\xi)d\tau$$

$$-\int_0^T \int_{\Omega(\xi)} \frac{\partial}{\partial \xi_i}\left(\nu(\xi,\tau)\frac{\partial C(\xi,\tau)}{\partial \xi_i}\right) G^*(\xi,x,\tau,t) d\Omega(\xi) d\tau$$

$$= -\int_0^T \int_{\Omega(\xi)} \frac{\partial}{\partial \xi_i}\left(\nu(\xi,\tau)\frac{\partial C(\xi,\tau)}{\partial \xi_i} G^*(\xi,x,\tau,t)\right) d\Omega(\xi) d\tau + \int_0^T \int_{\Omega(\xi)} \nu(\xi,\tau)\frac{\partial C(\xi,\tau)}{\partial \xi_i}\frac{\partial G^*(\xi,x,\tau,t)}{\partial \xi_i} d\Omega(\xi) d\tau$$

$$= -\int_0^T \int_{\Gamma(\xi)} \nu(\xi,\tau)\frac{\partial C(\xi,\tau)}{\partial n(\xi)} G^*(\xi,x,\tau,t) d\Gamma(\xi) d\tau + \int_0^T \int_{\Omega(\xi)} \frac{\partial}{\partial \xi_i}\left(C(\xi,\tau)\nu(\xi,\tau)\frac{\partial G^*(\xi,x,\tau,t)}{\partial \xi_i}\right) d\Omega(\xi) d\tau \quad (22c)$$

$$-\int_0^T \int_{\Omega(\xi)} C(\xi,\tau)\frac{\partial}{\partial \xi_i}\left(\nu(\xi,\tau)\frac{\partial G^*(\xi,x,\tau,t)}{\partial \xi_i}\right) d\Omega(\xi) d\tau$$

$$= -\int_0^T \int_{\Gamma(\xi)} \nu(\xi,\tau)\frac{\partial C(\xi,\tau)}{\partial n(\xi)} G^*(\xi,x,\tau,t) d\Gamma(\xi) d\tau$$

$$+ \int_0^T \int_{\Gamma(\xi)} C(\xi,\tau)\nu(\xi,\tau)\frac{\partial G^*(\xi,x,\tau,t)}{\partial n(\xi)} d\Gamma(\xi) d\tau - \int_0^T \int_{\Omega(\xi)} C(\xi,\tau)\frac{\partial}{\partial \xi_i}\left(\nu(\xi,\tau)\frac{\partial G^*(\xi,x,\tau,t)}{\partial \xi_i}\right) d\Omega(\xi) d\tau$$

Substituting Eq. (22) into Eq. (18), we have

$$0 = -\int_{\Omega(\xi)} C(\xi,0) G^*(\xi,x,0,t) d\Omega(\xi)$$

$$+ \int_0^T \int_{\Gamma(\xi)} \left[\begin{array}{l}\left(U_i(\xi,\tau)n_i(\xi,\tau)C(\xi,\tau) - \nu(\xi,\tau)\dfrac{\partial C(\xi,\tau)}{\partial n(\xi)}\right) G^*(\xi,x,\tau,t) \\ + C(\xi,\tau)\nu(\xi,\tau)\dfrac{\partial G^*(\xi,x,\tau,t)}{\partial n(\xi)}\end{array}\right] d\Gamma(\xi) d\tau \quad (23)$$

$$- \int_0^T \int_{\Omega(\xi)} \left[\begin{array}{l}C(\xi,\tau)\left(\dfrac{\partial G^*(\xi,x,\tau,t)}{\partial \tau} + \dfrac{\partial U_i(\xi,\tau) G^*(\xi,x,\tau,t)}{\partial \xi_i} + \dfrac{\partial}{\partial \xi_i}\left(\nu(\xi,\tau)\dfrac{\partial G^*(\xi,x,\tau,t)}{\partial \xi_i}\right)\right) \\ + \sigma(\xi,\tau) G^*(\xi,x,\tau,t)\end{array}\right] d\Omega(\xi) d\tau$$

Substituting Eqs. (16), (17) and (21b) into Eq. (23), we obtain

$$0 = -\int_{\Omega(\xi)} C_0(\xi) G(x,\xi,t,0) d\Omega(\xi)$$

$$+ \int_0^t \int_{\Gamma(\xi)} \left[\begin{array}{l}\left(U_i(\xi,\tau)n_i(\xi,\tau)C(\xi,\tau) - \nu(\xi,\tau)[b(\xi,\tau) - a(\xi,\tau)C(\xi,\tau)]\right) G(x,\xi,t,\tau) \\ + C(\xi,\tau)\nu(\xi,\tau)\dfrac{\partial G(x,\xi,t,\tau)}{\partial n(\xi)}\end{array}\right] d\Gamma(\xi) d\tau \quad , \quad (24)$$

$$- \int_0^t \int_{\Omega(\xi)} \left[\begin{array}{l}C(\xi,\tau)\left(\dfrac{\partial G(x,\xi,t,\tau)}{\partial \tau} + \dfrac{\partial U_i(\xi,\tau) G(x,\xi,t,\tau)}{\partial \xi_i} + \dfrac{\partial}{\partial \xi_i}\left(\nu(\xi,\tau)\dfrac{\partial G(x,\xi,t,\tau)}{\partial \xi_i}\right)\right) \\ + \sigma(\xi,\tau) G(x,\xi,t,\tau)\end{array}\right] d\Omega(\xi) d\tau$$

where the integral limit $T$ is replaced by $t$ because of Heaviside function $H(t-\tau)$ in $G$. Eq. (24) is a Volterra-Fredholm Integral Equation with unknown $C(x,t)$. Namely, $C(x,t)$, $(x,t) \in \Gamma(x) \times (0,T) + \Omega(x) \times (0,T)$ can be determined solving Eq. (24).

If the diffusion coefficient $\nu(x,t)$ is constant and the fundamental solution (FS) $G(x,\xi,t,\tau)$ given by Eq. (20) is used, we obtain

$$\varepsilon(x,t)C(x,t) = \int_{\Omega(\xi)} C_0(\xi) G(x,\xi,t,0) d\Omega(\xi)$$

$$- \int_0^t \int_{\Gamma(\xi)} \left[\begin{array}{l}\left(U_i(\xi,\tau)n_i(\xi,\tau)C(\xi,\tau) - \nu[b(\xi,\tau) - a(\xi,\tau)C(\xi,\tau)]\right) G(x,\xi,t,\tau) \\ + C(\xi,\tau)\nu\dfrac{\partial G(x,\xi,t,\tau)}{\partial n(\xi)}\end{array}\right] d\Gamma(\xi) d\tau \quad (25)$$

$$+ \int_0^t \int_{\Omega(\xi)} \left[C(\xi,\tau)\dfrac{\partial U_i(\xi,\tau) G(x,\xi,t,\tau)}{\partial \xi_i} + \sigma(\xi,\tau) G(x,\xi,t,\tau)\right] d\Omega(\xi) d\tau$$

where $\partial G/\partial \tau = -\partial G/\partial t$ and $\partial G/\partial \xi_i = -\partial G/\partial x_i$ are used, and $\varepsilon(x,t)$ is defined as

$$\varepsilon(x,t) = \begin{cases} 1 & \text{if } (x,t) \in \Omega \times (0,T) \\ 1/2 & \text{if } x \in \Omega \text{ and } t = 0 \text{ or } T, (x,t) \in \Gamma \times (0,T) \\ 0 & \text{otherwise} \end{cases} \quad (26)$$

Furthermore, if $U_i(x,t) = 0$, we have

$$\varepsilon(x)C(x,t) = \int_{\Omega(\xi)} C_0(\xi) G(x,\xi,t,0)\,d\Omega(\xi)$$
$$-\nu \int_0^t \int_{\Gamma(\xi)} \left[ -(b(\xi,\tau) - a(\xi,\tau)C(\xi,\tau))G(x,\xi,t,\tau) + C(\xi,\tau)\frac{\partial G(x,\xi,t,\tau)}{\partial n(\xi)} \right] d\Gamma(\xi)\,d\tau \quad (27)$$
$$+ \int_0^t \int_{\Omega(\xi)} \sigma(\xi,\tau) G(x,\xi,t,\tau)\,d\Omega(\xi)\,d\tau$$

When $(x,t) \in \Gamma \times (0,T)$, Eq. (27) becomes

$$\frac{1}{2}C(x,t) = \int_{\Omega(\xi)} C_0(\xi) G(x,\xi,t,0)\,d\Omega(\xi)$$
$$-\nu \int_0^t \int_{\Gamma(\xi)} \left[ -(b(\xi,\tau) - a(\xi,\tau)C(\xi,\tau))G(x,\xi,t,\tau) + C(\xi,\tau)\frac{\partial G(x,\xi,t,\tau)}{\partial n(\xi)} \right] d\Gamma(\xi)\,d\tau \quad (28)$$
$$+ \int_0^t \int_{\Omega(\xi)} \sigma(\xi,\tau) G(x,\xi,t,\tau)\,d\Omega(\xi)\,d\tau$$

Eq. (28) is a Volterra-Fredholm Integral Equation with unknown $C(x,t)$, $(x,t) \in \Gamma \times (0,T)$. Namely, $C(x,t)$, $(x,t) \in \Gamma \times (0,T)$ can be determined solving Eq. (28). $C(x,t)$, $(x,t) \in \Omega \times (0,T)$ can be obtained by Eq. (27).

## 4. Numerical results

### 4.1. 2D Space-Time Diffusion Problem

#### *4.1.1. Dirichlet Problem*

The boundary condition is given by
$$C(x,t) = c(x,t), \quad x = \pm L, \quad 0 < t < T. \quad (29)$$

The integral representation is given by

$$\varepsilon(x,t)C(x,t) = \int_{-L}^{+L} C_0(\xi) G(x,\xi,t,0)\,d\xi$$
$$- \int_0^t \left[ \begin{array}{l} +\nu \dfrac{\partial C(-L,t)}{\partial \xi} G(x,-L,t,\tau) + \nu c(-L,\tau) \dfrac{\partial G(x,-L,t,\tau)}{\partial x} \\ -\nu \dfrac{\partial C(+L,t)}{\partial \xi} G(x,+L,t,\tau) - \nu c(+L,\tau) \dfrac{\partial G(x,+L,t,\tau)}{\partial x} \end{array} \right] d\tau \quad (30)$$
$$+ \int_0^t \int_{-L}^{+L} \sigma(\xi,\tau) G(x,\xi,t,\tau)\,d\xi\,d\tau,$$
$$-L \leq x \leq +L, \quad 0 < t < T.$$

The integral equation is given by

$$\int_0^t \left[ +\nu C_\xi(-L,t) G(-L,-L,t,\tau) - \nu C_\xi(+L,t) G(-L,+L,t,\tau) \right] d\tau$$
$$= -c(-L,t)/2 + \int_{-L}^{+L} C_0(\xi) G(-L,\xi,t,0)\,d\xi + \int_0^t \int_{-L}^{+L} \sigma(\xi,\tau) G(-L,\xi,t,\tau)\,d\xi\,d\tau \quad (31a)$$
$$- \int_0^t \left[ +\nu c(-L,\tau) G_x(-L,-L,t,\tau) - \nu c(+L,\tau) G_x(-L,+L,t,\tau) \right] d\tau,$$

$$\int_0^t \left[ +\nu C_\xi(-L,t) G(+L,-L,t,\tau) - \nu C_\xi(+L,t) G(+L,+L,t,\tau) \right] d\tau$$
$$= -c(+L,t)/2 + \int_{-L}^{+L} C_0(\xi) G(+L,\xi,t,0)\,d\xi + \int_0^t \int_{-L}^{+L} \sigma(\xi,\tau) G(+L,\xi,t,\tau)\,d\xi\,d\tau \quad (31b)$$
$$- \int_0^t \left[ +\nu c(-L,\tau) G_x(+L,-L,t,\tau) - \nu c(+L,\tau) G_x(+L,+L,t,\tau) \right] d\tau.$$

Eq. (31) is a Simultaneous Volterra Integral Equation with unknown $C_\xi(\pm L,\tau)$.

The exact solution of Dirichlet problem is given in Appendix C.1.
Numerical parameters are as follows:

$$C_0(x) = \exp\left([x/(L/8)]^2\right), \quad \sigma(x,t) = 0, \quad c(\pm L, t) = 0,$$
$$v = 0.05, \quad L = 1, \quad dx = L/41 = 0.04878, \quad dt = 0.0625. \tag{32}$$

The numerical results is shown in Fig. 4. The numerical result (C) shows a satisfactory agreement with the exact solution (Cex).

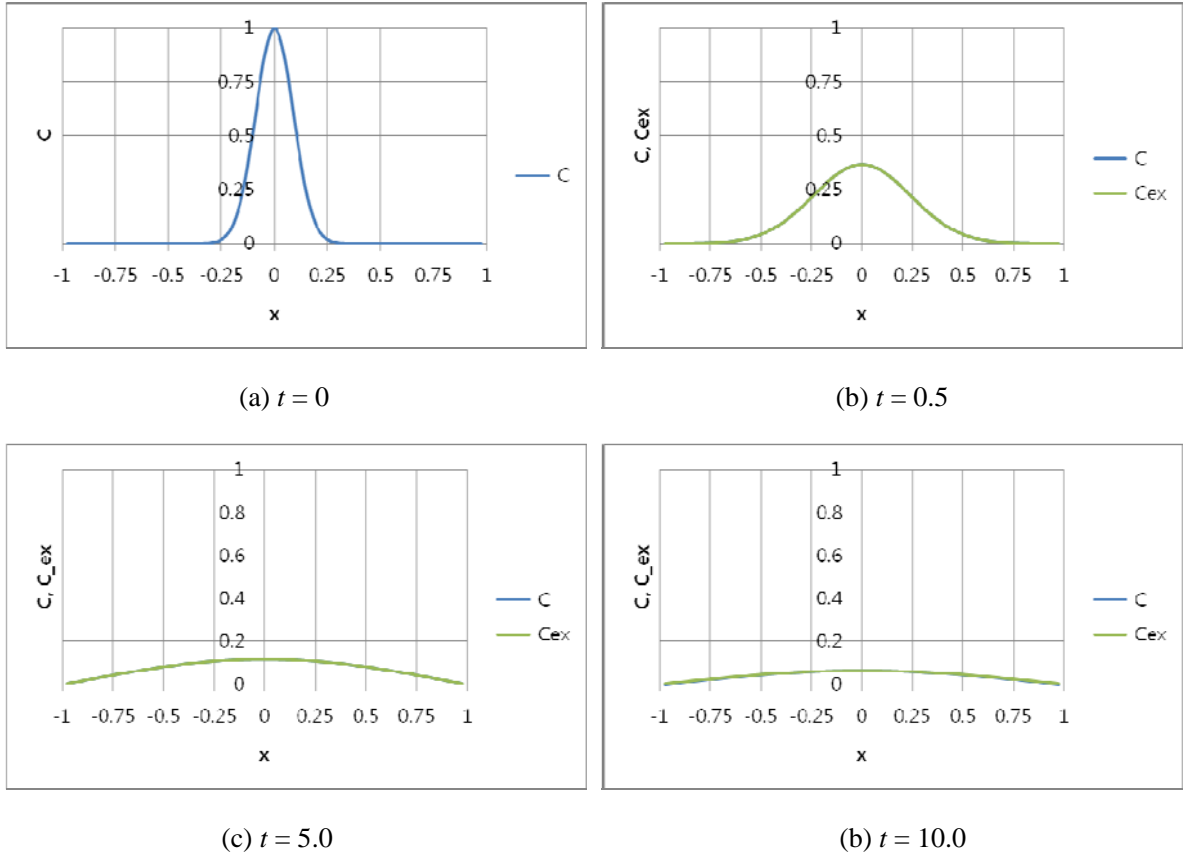

(a) $t = 0$      (b) $t = 0.5$

(c) $t = 5.0$      (b) $t = 10.0$

Fig. 4. Numerical results of Dirichlet problem

### 4.1.2. Neumann Problem

The boundary condition is given by
$$\frac{\partial C(x,t)}{\partial x} = b(x,t), \quad x = \pm L, \quad 0 < t < T. \tag{32}$$

The integral representation is given by

$$\varepsilon(x,t)C(x,t) = \int_{-L}^{+L} C_0(\xi) G(x,\xi,t,0)\, d\xi$$
$$- \int_0^t \left[ \begin{array}{l} + v\, b(-L,\tau)\, G(x,-L,t,\tau) + v\, C(-L,\tau)\dfrac{\partial G(x,-L,t,\tau)}{\partial x} \\ - v\, b(+L,\tau)\, G(x,+L,t,\tau) - v\, C(+L,\tau)\dfrac{\partial G(x,+L,t,\tau)}{\partial x} \end{array} \right] d\tau \tag{33}$$
$$+ \int_0^t \int_{-L}^{+L} \sigma(\xi,\tau)\, G(x,\xi,t,\tau)\, d\xi\, d\tau,$$
$$-L \le x \le +L, \quad 0 < t < T.$$

The integral equation is given by

$$C(-L,t)/2 - \int_0^t [+\nu C(-L,\tau)G_x(-L,-L,t,\tau) - \nu C(+L,\tau)G_x(-L,+L,t,\tau)]d\tau$$
$$= +\int_{-L}^{+L} C_0(\xi)G(-L,\xi,t,0)d\xi + \int_0^t \int_{-L}^{+L} \sigma(\xi,\tau)G(-L,\xi,t,\tau)d\xi d\tau \tag{34a}$$
$$- \int_0^t [+\nu b(-L,\tau)G(-L,-L,t,\tau) - \nu b(+L,\tau)G(-L,+L,t,\tau)]d\tau,$$

$$C(+L,t)/2 - \int_0^t [+\nu C(-L,\tau)G_x(+L,-L,t,\tau) - \nu C(+L,\tau)G_x(+L,+L,t,\tau)]d\tau$$
$$= +\int_{-L}^{+L} C_0(\xi)G(+L,\xi,t,0)d\xi + \int_0^t \int_{-L}^{+L} \sigma(\xi,\tau)G(+L,\xi,t,\tau)d\xi d\tau \tag{34b}$$
$$- \int_0^t [+\nu b(-L,\tau)G(+L,-L,t,\tau) - \nu b(+L,\tau)G(+L,+L,t,\tau)]d\tau.$$

Eq. (34) is a Simultaneous Volterra Integral Equation with unknown $C(\pm L,\tau)$.

The exact solution of Neumann problem is given in Appendix C.2.

Numerical parameters are as follows:
$$C_0(x) = \exp([x/(L/8)]^2), \quad \sigma(x,t) = 0, \quad b(\pm L,t) = 0,$$
$$\nu = 0.05, \quad L = 1, \quad M = 161, \quad dx = 2L/M = 0.01242, \quad dt = 0.005. \tag{35}$$

The numerical results is shown in Fig. 5. The numerical result (C) shows a satisfactory agreement with the exact solution (Cex). The time increment $dt$ should be taken smaller for Neumann problem than for Dirichlet problem. When time $t$ is small, $M$ should be big. On the contrary, when $t$ is big, small $M$ stabilizes the calculation.

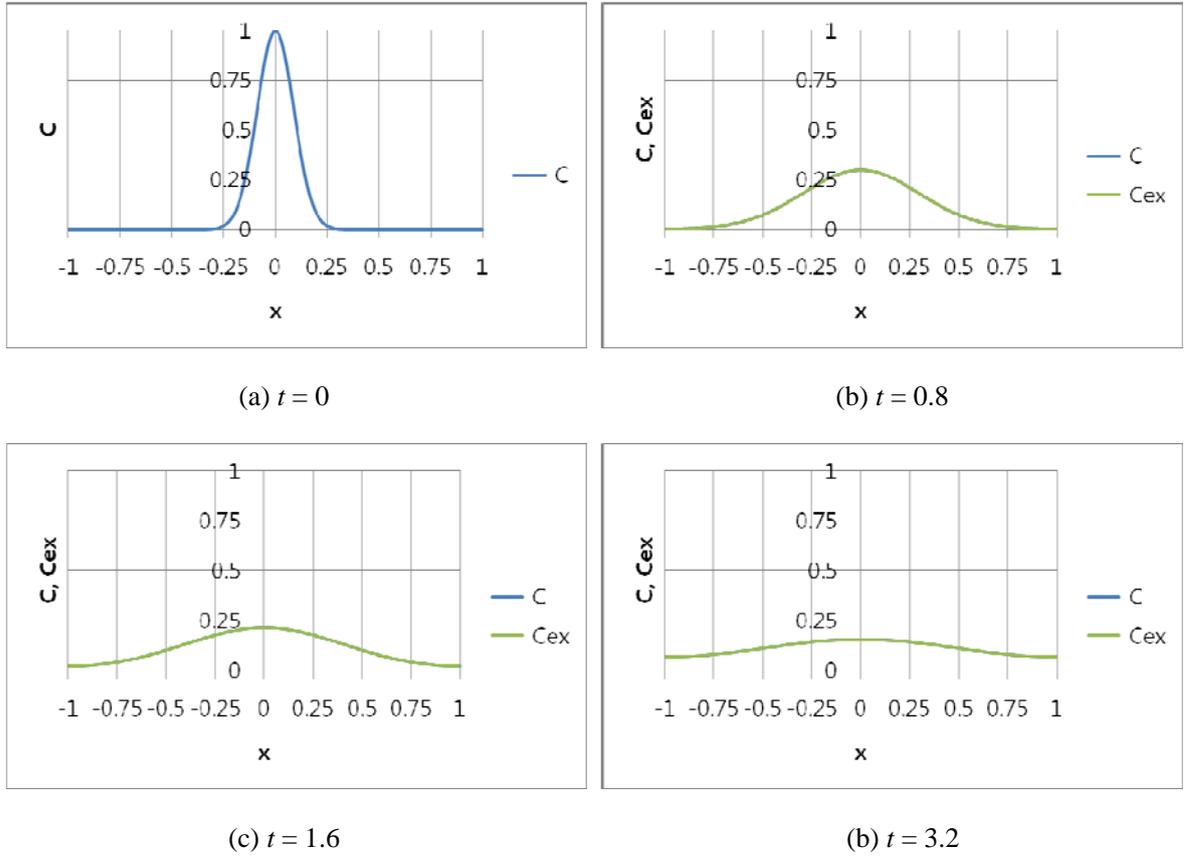

(a) $t = 0$  (b) $t = 0.8$

(c) $t = 1.6$  (b) $t = 3.2$

Fig. 5. Numerical results of Neumann problem

## 5. Conclusions

(1) We have discussed on the Generalized Integral Representation Method (GIRM) for Space-Time-Separated

Method (STSM) and Space-Time-Unified Method (STUM).
(2) STSM and STUM give explicit and implicit time evolutions, respectively.
(3) The algorithm of STSM is much simpler than STUM. However, the implicit time evolution could give us much more efficient computation.
(4) We have conducted numerical calculations using STUM for Dirichlet and Neumann problems in 2D space-time. We used a Traditional Fundamental Solution (TFS). The results seem very satisfactory.
(5) However, in case of Neumann problem, the time increment becomes much smaller than in case of Dirichlet problem.
(6) We have studied the reason. This property may be explained by the nature of the Traditional Fundamental Solution (TFS) used for the numerical calculations.

**References**

1. J.C. Wu, J,F. "Thompson, "Numerical solutions of time-dependent incompressible Navier-Stokes equations using an integro-differential formulations", *Computers & Fluids*, (1973), 1, pp. 197-215.
2. J.S. Uhlman, An integral equation formulation of the equations of motion of an incompressible fluid, NUWC-NPT Technical Report 10,086 15 July (1992).
3. Yang Hongtao, On the Convergence of Boundary Element Methods for Initial-Neumann Problems for the Heat Equation, Mathematics of Computation, Vol. 68, No. 226, April 1999, pp. 547-557.
4. H. Isshik, S. Nagata, Y. Imai, Solution of Viscous Flow around a Circular Cylinder by a New Integral Representation Method (NIRM), Asian Journal of Engineering and Technology (AJET), Vol. 02, Issue o2, April 2014, pp. 60-82.
5. H. Isshiki, From Integral Representation Method (IRM) to Generalized Integral Representation Method (GIRM), Applied and Computational Mathematics, Special Issue: Integral Representation Method and Its Generalization, Vol. 4, No. 3-1, 2015, pp. 1-14.
6. G. Tsedendorj, H. Isshiki, Numerical Solution of Two-Dimensional Advection-Diffusion Using Generalized Integral Representation Method, International Journal of Computational Methods, Vol. 14, No. 1, 2017, pp. 1750028-1-13

**Appendix A. A Fundamental Solution for Differential Operator $\Delta$ Used for Steady Diffusion Problem**

A fundamental solution $G(x_1, x_2, t)$ having the singularity at $x_1 = x = t = 0$ is defined as

$$\frac{\partial^2 G(x_1, x_2)}{\partial x_1^2} + \frac{\partial^2 G(x_1, x_2)}{\partial x_2^2} = \delta(x_1)\delta(x_2), \quad (x_1, x_2) \in R^2 \tag{A1}$$

Eq. (A1) can be written in polar coordinates:

$$x_1 = r\cos\theta, \quad x_2 = r\sin\theta \tag{A2}$$

as

$$\frac{1}{r}\frac{\partial}{\partial r}\left(r\frac{\partial G}{\partial r}\right) + \frac{1}{r^2}\frac{\partial^2 G}{\partial \theta^2} = 0 \tag{A3}$$

A singular solution having the lowest singularity:

$$G = \frac{1}{2\pi}\ln r \tag{A4}$$

is the fundamental, since $G$ satisfies Eq. (A1) because

$$\int_0^R \int_0^{2\pi} \left(\frac{\partial^2 G(x_1, x_2)}{\partial x_1^2} + \frac{\partial^2 G(x_1, x_2)}{\partial x_2^2}\right) r\, dr\, d\theta = \int_0^{2\pi}\left[\frac{\partial G}{\partial r}\right]_{r=R} R\, d\theta = 1 \tag{A5}$$

**Appendix B. A Fundamental Solution for Differential Operator $\partial/\partial t - \Delta$ Used for Unsteady Diffusion Problem**

A fundamental solution $G(x_1, x_2, t)$ having the singularity at $x_1 = x = t = 0$ is defined as

$$\frac{\partial G(x_1, x_2, t)}{\partial t} - \nu\left(\frac{\partial^2 G(x_1, x_2, t)}{\partial x_1^2} + \frac{\partial^2 G(x_1, x_2, t)}{\partial x_2^2}\right) = \delta(x_1)\delta(x_2)\delta(t) \tag{B1}$$

Fourier transforms with respects to $x_1$ and $x_2$ are given by

$$\hat{G}(k_1, k_2, t) = \frac{1}{2\pi} \int_{-\infty}^{+\infty}\int_{-\infty}^{+\infty} G(x_1, x_2, t) e^{-ik_1 x_1 - ik_2 x_2} dx_1 dx_2 \tag{B2a}$$

$$G(x_1, x_2, t) = \frac{1}{2\pi} \int_{-\infty}^{+\infty}\int_{-\infty}^{+\infty} \hat{G}(k_1, k_2, t) e^{ik_1 x_1 + ik_2 x_2} dk_1 dk_2 \tag{B2b}$$

Operating $\frac{1}{2\pi} \int_{-\infty}^{+\infty}\int_{-\infty}^{+\infty} dx_1 dx_2 e^{-ik_1 x_1 - ik_2 x_2} \cdot$ to the both sides of Eq. (B1), we obtain

$$\frac{\partial \hat{G}(k_1, k_2, t)}{\partial t} - \nu(k_1^2 + k_2^2)\hat{G}(k_1, k_2, t) = \frac{1}{2\pi}\delta(t) \tag{B3}$$

From Eq. (B3), we obtain

$$\hat{G}(k_1, k_2, t) = \frac{1}{2\pi} H(t) \exp\left(-\nu(k_1^2 + k_2^2)t\right) \tag{B4}$$

Using Eq. (B2b), we have

$$\begin{aligned} G(x_1, x_2, t) &= \frac{1}{4\pi^2} H(t) \int_{-\infty}^{+\infty}\int_{-\infty}^{+\infty} \exp\left(-\nu(k_1^2 + k_2^2)t\right) e^{ik_1 x_1 + ik x_2} dk_1 dk_2 \\ &= \frac{1}{4\pi^2} H(t) \int_{-\infty}^{+\infty} \exp\left(-\nu k_1^2 t + ik_1 x_1\right) dk_1 \int_{-\infty}^{+\infty} \exp\left(-\nu k_2^2 t + ik_2 x_2\right) dk_2 \end{aligned} \tag{B5}$$

Since an integral formula gives us

$$\int_{-\infty}^{+\infty} \exp\left(-\nu k_1^2 t + ik_1 x_1\right) dk_1 = \int_{-\infty}^{+\infty} \exp\left(-\nu t\left(k_1^2 + ik_1 \frac{x_1}{\nu t}\right)\right) dk_1 = \sqrt{\frac{\pi}{\nu t}} \exp\left(-\frac{x_1^2}{4\nu t}\right) \tag{B6}$$

We derive

$$G(x_1, x_2, t) = \frac{1}{4\pi\nu t} H(t) \exp\left(-\frac{x_1^2 + x_2^2}{4\nu t}\right) \tag{B7}$$

## Appendix C. Exact Solution of Diffusion Problem

### C.1. Dirichlet Problem

The initial-boundary problem of Drichlet problem is given by

$$\frac{\partial C(x,t)}{\partial t} = \nu \frac{\partial^2 C(x,t)}{\partial x^2}, \quad -L < x < L, \ 0 < t < T, \tag{C1}$$

$$C(\pm L, t) = g_{\pm L}(t), \quad 0 < t < T, \tag{C2}$$

$$C(x, 0) = C_0(x), \quad -L < x < +L. \tag{C3}$$

We assume

$$C(x,t) = \frac{g_{-L}(t)}{2L}(L-x) + \frac{g_{+L}(t)}{2L}(L+x) + \sum_{n=1}^{N} A_n(t) \sin\frac{n\pi}{2L}(x+L), \tag{C4}$$

$$L \pm x = \sum_{n=1}^{N} a^{\pm}_n \sin\frac{n\pi}{2L}(x+L), \tag{C5}$$

where

$$a^{\pm}_m = \frac{2L}{m\pi}\left[-(-1)^m + 1\right](1 \pm 1). \tag{C6}$$

Substituting Eq. (C4) and (C5) into Eq. (C1), we have

$$\begin{aligned} &\frac{1}{2L}\frac{dg_{-L}(t)}{dt}\sum_{n=1}^{N} a^{-}_n \sin\frac{n\pi}{2L}(x+L) + \frac{1}{2L}\frac{dg_{+L}(t)}{dt}\sum_{n=1}^{N} a^{+}_n \sin\frac{n\pi}{2L}(x+L) \\ &+ \sum_{n=1}^{N}\left(\frac{dA_n(t)}{dt} + \nu\left(\frac{n\pi}{2L}\right)^2 A_n(t)\right)\sin\frac{n\pi}{2L}(x+L) = 0. \end{aligned} \tag{C7}$$

From Eq. (C7), we obtain

$$\frac{dA_n(t)}{dt} + \nu\left(\frac{n\pi}{2L}\right)^2 A_n(t) = -\frac{1}{2L}\frac{dg_{-L}(t)}{dt}a^{-}_n - \frac{1}{2L}\frac{dg_{+L}(t)}{dt}a^{+}_n. \tag{C8}$$

From Eqs. (C3), (C4) and (C5), we have

$$\frac{g_{-L}(0)}{2L}\sum_{n=1}^{N} a^{-}_n \sin\frac{n\pi}{2L}(x+L) + \frac{g_{+L}(0)}{2L}\sum_{n=1}^{N} a^{+}_n \sin\frac{n\pi}{2L}(x+L) + \sum_{n=1}^{N} A_n(0)\sin\frac{n\pi}{2L}(x+L) = C_0(x). \tag{C9}$$

Let $C_0(x)$ be

$$C_0(x) = \sum_{n=1}^{N} c_n \sin \frac{n\pi}{2L}(x+L), \tag{C10}$$

where

$$c_m = \frac{1}{L}\int_{-L}^{+L} C_0(x) \sin \frac{m\pi}{2L}(x+L)\, dx. \tag{C11}$$

From Eqs. (C9) and (C10), the initial condition of Eq. (C8) is given by

$$A_n(0) = c_n - \frac{g_{-L}(0)}{2L} a^-_n - \frac{g_{+L}(0)}{2L} a^+_n. \tag{C12}$$

Then, the solution of Eqs. (C8) and (C12) is given by

$$A_n(t) = \left(c_n - \frac{g_{-L}(0)}{2L} a^-_n - \frac{g_{+L}(0)}{2L} a^+_n\right) G_n(t) - \int_0^t G_n(t-\tau)\left(\frac{1}{2L}\frac{dg_{-L}(\tau)}{d\tau} a^-_n + \frac{1}{2L}\frac{dg_{+L}(\tau)}{dt} a^+_n\right) d\tau, \tag{C13}$$

where

$$G_n(t) = \exp\left(-\nu\left(\frac{n\pi}{2L}\right)^2 t\right) H(t). \tag{C14}$$

## C.2. Neumann Problem

The initial-boundary problem of Newmann problem is given by

$$\frac{\partial C(x,t)}{\partial t} = \nu \frac{\partial^2 C(x,t)}{\partial x^2}, \quad -L < x < L,\ 0 < t < T \tag{C15}$$

$$\frac{\partial C(\pm L, t)}{\partial x} = \pm f_{\pm L}(t), \quad 0 < t < T \tag{C16}$$

$$C(x,0) = C_0(x), \quad -L < x < +L \tag{C17}$$

We assume

$$C(x,t) = -\frac{f_{-L}(t)}{4L}(L-x)^2 + \frac{f_{+L}(t)}{4L}(L+x)^2 + \sum_{n=0}^{N} B_n(t) \cos \frac{n\pi}{2L}(x+L), \tag{C18}$$

$$(L \pm x)^2 = \sum_{n=0}^{N} b^\pm_n \cos \frac{n\pi}{2L}(x+L), \tag{C19}$$

where

$$b^\pm_m = 2\left[-(1\mp 1)\left((-1)^m - 1\right) + 2(-1)^m\right]\left(\frac{2L}{m\pi}\right)^2 \frac{1}{\varepsilon_{mm}} \tag{C20}$$

and $\varepsilon_{nm} = 2$ when $n = m = 0$, 1 otherwise.

Substituting Eq. (C18) and (C19) into Eq. (C15), we have

$$-\frac{1}{4L}\frac{df_{-L}(t)}{dt}\left(\sum_{n=0}^{N} b^-_n \cos\frac{n\pi}{2L}(x+L) - 2\nu\right) + \frac{1}{4L}\frac{df_{+L}(t)}{dt}\left(\sum_{n=0}^{N} b^+_n \cos\frac{n\pi}{2L}(x+L) - 2\nu\right)$$

$$+ \sum_{n=0}^{N}\left(\frac{dB_n(t)}{dt} + \nu\left(\frac{n\pi}{2L}\right)^2 B_n(t)\right)\cos\frac{n\pi}{2L}(x+L) = 0. \tag{C21}$$

From Eq. (C21), we obtain

$$\frac{dB_n(t)}{dt} + \nu\left(\frac{n\pi}{2L}\right)^2 B_n(t) = \frac{1}{4L}\frac{df_{-L}(t)}{dt}\left(b^-_n - 2\nu\delta_{n0}\right) - \frac{1}{4L}\frac{df_{+L}(t)}{dt}\left(b^+_n - 2\nu\delta_{n0}\right). \tag{C22}$$

From Eqs. (C17), (C18) and (C19), we have

$$C(x,0) = -\frac{f_{-L}(0)}{4L}\sum_{n=0}^{N} b^-_n \cos\frac{n\pi}{2L}(x+L) + \frac{f_{+L}(0)}{4L}\sum_{n=0}^{N} b^+_n \cos\frac{n\pi}{2L}(x+L) + \sum_{n=0}^{N} B_n(0) \cos\frac{n\pi}{2L}(x+L) = C_0(x). \tag{C23}$$

Let $C_0(x)$ be

$$C_0(x) = \sum_{n=0}^{N} c'_n \cos\frac{n\pi}{2L}(x+L), \tag{C24}$$

where

$$c'_m = \frac{1}{L\varepsilon_{mm}} \int_{-L}^{+L} C_0(x) \cos\frac{m\pi}{2L}(x+L)\, dx. \tag{C25}$$

From Eqs. (C23) and (C24), the initial condition of Eq. (C22) is given by

$$B_n(0) = c'_n + \frac{f_{-L}(0)}{4L} b^-_n - \frac{f_{+L}(0)}{4L} b^+_n. \tag{C26}$$

Then, the solution of Eqs. (C22) and (C26) is given by

$$\begin{aligned}B_n(t) =& \left( c'_n + \frac{f_{-L}(0)}{4L} b^-_n - \frac{f_{+L}(0)}{4L} b^+_n \right) G_n(t) \\ & - \int_0^t G_n(t-\tau) \left( -\frac{1}{4L} \frac{df_{-L}(\tau)}{d\tau} (b^-_n - 2\nu\delta_{n0}) + \frac{1}{4L} \frac{df_{+L}(\tau)}{dt} (b^+_n - 2\nu\delta_{n0}) \right) d\tau\end{aligned}, \tag{C27}$$

where

$$G_n(t) = \exp\left( -\nu \left( \frac{n\pi}{2L} \right)^2 t \right) H(t). \tag{C28}$$